\documentclass{article}

\usepackage{PRIMEarxiv}

\usepackage[utf8]{inputenc} 
\usepackage[T1]{fontenc}    
\usepackage{hyperref}       
\usepackage{url}            
\usepackage{booktabs}       
\usepackage{amsfonts}       
\usepackage{nicefrac}       
\usepackage{microtype}      
\usepackage{lipsum}
\usepackage{fancyhdr}       
\usepackage{graphicx}       
\graphicspath{{media/}}     

\usepackage{comment}
\usepackage{graphicx}
\usepackage{amsmath}
\usepackage{amsfonts}
\usepackage{listings}
\usepackage[version=4]{mhchem}
\usepackage{siunitx}
\usepackage{longtable,tabularx}

\newcommand{\br}{\textbf{\emph{r}}}

\newcommand{\bv}{\textbf{\emph{v}}}

\newcommand{\RR}{\mathbb{R}}

\usepackage{listings}
\usepackage{xcolor}
\definecolor{maccolor}{rgb}{0.3,0.3,0.8}
\lstdefinelanguage{Macaulay2}
{
basicstyle={\ttfamily},
keywordstyle={\color{maccolor!80!black}},
commentstyle={\color{gray}},
stringstyle={\color{red!40!black}},
rulecolor=\color{maccolor},
basewidth={1.2ex}, 
sensitive=false,
morecomment=[l]{--},
morecomment=[s]{-*}{*-},
morestring=[b]",
escapechar={`},
escapebegin={\rmfamily},
morekeywords={about,abs,AbstractToricVarieties,accumulate,Acknowledgement,acos,acosh,acot,addCancelTask,addDependencyTask,addEndFunction,addHook,AdditionalPaths,addStartFunction,addStartTask,Adjacent,adjoint,AdjointIdeal,AffineVariety,AfterEval,AfterNoPrint,AfterPrint,agm,AInfinity,alarm,AlgebraicSplines,Algorithm,Alignment,all,AllCodimensions,allowableThreads,ambient,analyticSpread,Analyzer,AnalyzeSheafOnP1,ancestor,ancestors,ANCHOR,and,andP,AngleBarList,ann,annihilator,antipode,any,append,applicationDirectory,applicationDirectorySuffix,apply,applyKeys,applyPairs,applyTable,applyValues,apropos,argument,Array,arXiv,Ascending,ascii,asin,asinh,ass,assert,associatedGradedRing,associatedPrimes,AssociativeAlgebras,AssociativeExpression,atan,atan2,atEndOfFile,Authors,autoload,AuxiliaryFiles,backtrace,Bag,Bareiss,baseFilename,BaseFunction,baseName,baseRing,baseRings,BaseRow,BasicList,basis,BasisElementLimit,Bayer,BeforePrint,beginDocumentation,BeginningMacaulay2,Benchmark,benchmark,Bertini,BesselJ,BesselY,betti,BettiCharacters,BettiTally,between,BGG,BIBasis,Binary,BinaryOperation,Binomial,binomial,BinomialEdgeIdeals,Binomials,BKZ,BlockMatrix,BLOCKQUOTE,BODY,Body,BoijSoederberg,BOLD,Book3264Examples,Boolean,BooleanGB,borel,Boxes,BR,break,Browse,Bruns,cache,CacheExampleOutput,CacheFunction,CacheTable,cacheValue,CallLimit,cancelTask,capture,catch,Caveat,CC,CDATA,ceiling,Center,centerString,Certification,ChainComplex,chainComplex,ChainComplexExtras,ChainComplexMap,ChainComplexOperations,ChangeMatrix,char,CharacteristicClasses,characters,charAnalyzer,check,CheckDocumentation,chi,Chordal,class,Classic,clean,clearAll,clearEcho,clearOutput,close,closeIn,closeOut,ClosestFit,CODE,code,codim,CodimensionLimit,coefficient,CoefficientRing,coefficientRing,coefficients,Cofactor,CohenEngine,CohenTopLevel,CoherentSheaf,CohomCalg,cohomology,coimage,CoincidentRootLoci,coker,cokernel,collectGarbage,columnAdd,columnate,columnMult,columnPermute,columnRankProfile,columnSwap,combine,Command,commandInterpreter,commandLine,COMMENT,commonest,commonRing,comodule,CompactMatrix,compactMatrixForm,CompiledFunction,CompiledFunctionBody,CompiledFunctionClosure,Complement,complement,complete,CompleteIntersection,CompleteIntersectionResolutions,Complexes,ComplexField,components,compose,compositions,compress,concatenate,conductor,ConductorElement,cone,Configuration,ConformalBlocks,conjugate,connectionCount,Consequences,Constant,Constants,constParser,content,continue,contract,Contributors,ConvexInterface,conwayPolynomial,ConwayPolynomials,copy,copyDirectory,copyFile,copyright,Core,CorrespondenceScrolls,cos,cosh,cot,CotangentSchubert,cotangentSheaf,coth,cover,coverMap,cpuTime,createTask,Cremona,csc,csch,current,currentColumnNumber,currentDirectory,currentFileDirectory,currentFileName,currentLayout,currentLineNumber,currentPackage,currentString,currentTime,Cyclotomic,Database,Date,DD,dd,deadParser,debug,debugError,DebuggingMode,debuggingMode,debugLevel,DecomposableSparseSystems,Decompose,decompose,deepSplice,Default,default,defaultPrecision,Degree,degree,degreeLength,DegreeLift,DegreeLimit,DegreeMap,DegreeOrder,DegreeRank,Degrees,degrees,degreesMonoid,degreesRing,delete,demark,denominator,Dense,Density,Depth,depth,Descending,Descent,Describe,describe,Description,det,determinant,DeterminantalRepresentations,DGAlgebras,diagonalMatrix,diameter,Dictionary,dictionary,dictionaryPath,diff,DiffAlg,difference,dim,directSum,disassemble,discriminant,dismiss,Dispatch,distinguished,DIV,Divide,divideByVariable,DivideConquer,DividedPowers,Divisor,DL,Dmodules,do,doc,docExample,docTemplate,document,DocumentTag,Down,drop,DT,dual,eagonNorthcott,EagonResolution,echoOff,echoOn,EdgeIdeals,edit,EigenSolver,eigenvalues,eigenvectors,eint,EisenbudHunekeVasconcelos,elapsedTime,elapsedTiming,elements,Eliminate,eliminate,Elimination,EliminationMatrices,EllipticCurves,EllipticIntegrals,else,EM,Email,End,end,endl,endPackage,Engine,engineDebugLevel,EngineRing,EngineTests,entries,EnumerationCurves,environment,Equation,EquivariantGB,erase,erf,erfc,error,errorDepth,euler,EulerConstant,eulers,even,EXAMPLE,ExampleFiles,ExampleItem,examples,ExampleSystems,Exclude,exec,exit,exp,expectedReesIdeal,expm1,exponents,export,exportFrom,exportMutable,Expression,expression,Ext,extend,ExteriorIdeals,ExteriorModules,exteriorPower,Factor,factor,false,Fano,FastMinors,FastNonminimal,FGLM,File,fileDictionaries,fileExecutable,fileExists,fileExitHooks,fileLength,fileMode,FileName,FilePosition,fileReadable,fileTime,fileWritable,fillMatrix,findFiles,findHeft,FindOne,findProgram,findSynonyms,FiniteFittingIdeals,First,first,firstkey,FirstPackage,fittingIdeal,flagLookup,FlatMonoid,flatten,flattenRing,Flexible,flip,floor,flush,fold,FollowLinks,for,forceGB,fork,FormalGroupLaws,Format,format,formation,FourierMotzkin,FourTiTwo,fpLLL,frac,fraction,FractionField,frames,FrobeniusThresholds,from,fromDividedPowers,fromDual,Function,FunctionApplication,FunctionBody,functionBody,FunctionClosure,FunctionFieldDesingularization,fusePairs,futureParser,GaloisField,Gamma,gb,GBDegrees,gbRemove,gbSnapshot,gbTrace,gcd,gcdCoefficients,gcdLLL,GCstats,genera,GeneralOrderedMonoid,GenerateAssertions,generateAssertions,generator,generators,Generic,GenericInitialIdeal,genericMatrix,genericSkewMatrix,genericSymmetricMatrix,gens,genus,get,getc,getChangeMatrix,getenv,getGlobalSymbol,getNetFile,getNonUnit,getPrimeWithRootOfUnity,getSymbol,getWWW,GF,gfanInterface,Givens,GKMVarieties,GLex,Global,global,globalAssign,globalAssignFunction,GlobalAssignHook,globalAssignment,globalAssignmentHooks,GlobalDictionary,GlobalHookStore,globalReleaseFunction,GlobalReleaseHook,Gorenstein,GradedLieAlgebras,GradedModule,gradedModule,GradedModuleMap,gradedModuleMap,gramm,GraphicalModels,GraphicalModelsMLE,Graphics,graphIdeal,graphRing,Graphs,Grassmannian,GRevLex,GroebnerBasis,groebnerBasis,GroebnerBasisOptions,GroebnerStrata,GroebnerWalk,groupID,GroupLex,GroupRevLex,GTZ,Hadamard,handleInterrupts,HardDegreeLimit,hash,HashTable,hashTable,HEAD,HEADER1,HEADER2,HEADER3,HEADER4,HEADER5,HEADER6,HeaderType,Heading,Headline,Heft,heft,Height,height,help,Hermite,hermite,Hermitian,HH,hh,HigherCIOperators,HighestWeights,Hilbert,hilbertFunction,hilbertPolynomial,hilbertSeries,HodgeIntegrals,hold,Holder,Hom,homeDirectory,HomePage,Homogeneous,Homogeneous2,homogenize,homology,homomorphism,HomotopyLieAlgebra,hooks,horizontalJoin,HorizontalSpace,HR,HREF,HTML,html,httpHeaders,Hybrid,HyperplaneArrangements,Hypertext,hypertext,HypertextContainer,HypertextParagraph,icFracP,icFractions,icMap,icPIdeal,id,Ideal,ideal,idealizer,identity,if,IgnoreExampleErrors,ii,image,imaginaryPart,IMG,ImmutableType,importFrom,in,incomparable,Increment,independentSets,indeterminate,IndeterminateNumber,Index,index,indexComponents,IndexedVariable,IndexedVariableTable,indices,inducedMap,inducesWellDefinedMap,InexactField,InexactFieldFamily,InexactNumber,InfiniteNumber,infinity,info,InfoDirSection,infoHelp,Inhomogeneous,input,Inputs,insert,installAssignmentMethod,installedPackages,installHilbertFunction,installMethod,installMinprimes,installPackage,InstallPrefix,instance,instances,IntegralClosure,integralClosure,integrate,IntermediateMarkUpType,interpreterDepth,intersect,intersectInP,Intersection,intersection,interval,InvariantRing,inverse,InverseMethod,inversePermutation,Inverses,inverseSystem,InverseSystems,Invertible,InvolutiveBases,irreducibleCharacteristicSeries,irreducibleDecomposition,isAffineRing,isANumber,isBorel,isCanceled,isCommutative,isConstant,isDirectory,isDirectSum,isEmpty,isField,isFinite,isFinitePrimeField,isFreeModule,isGlobalSymbol,isHomogeneous,isIdeal,isInfinite,isInjective,isInputFile,isIsomorphism,isLinearType,isListener,isLLL,isMember,isModule,isMonomialIdeal,isNormal,isOpen,isOutputFile,isPolynomialRing,isPrimary,isPrime,isPrimitive,isPseudoprime,isQuotientModule,isQuotientOf,isQuotientRing,isReady,isReal,isReduction,isRegularFile,isRing,isSkewCommutative,isSorted,isSquareFree,isStandardGradedPolynomialRing,isSubmodule,isSubquotient,isSubset,isSupportedInZeroLocus,isSurjective,isTable,isUnit,isWellDefined,isWeylAlgebra,ITALIC,Iterate,Jacobian,jacobian,jacobianDual,Jets,Join,join,Jupyter,K3Carpets,K3Surfaces,Keep,KeepFiles,KeepZeroes,ker,kernel,kernelLLL,kernelOfLocalization,Key,keys,Keyword,Keywords,kill,koszul,Kronecker,KustinMiller,LABEL,last,lastMatch,LATER,LatticePolytopes,Layout,lcm,leadCoefficient,leadComponent,leadMonomial,leadTerm,Left,left,length,LengthLimit,letterParser,Lex,LexIdeals,LI,Licenses,LieTypes,lift,liftable,Limit,limitFiles,limitProcesses,Linear,LinearAlgebra,LinearTruncations,lineNumber,lines,LINK,linkFile,List,list,listForm,listLocalSymbols,listSymbols,listUserSymbols,LITERAL,LLL,LLLBases,lngamma,load,loadDepth,LoadDocumentation,loadedFiles,loadedPackages,loadPackage,Local,local,localDictionaries,LocalDictionary,localize,LocalRings,locate,log,log1p,LongPolynomial,lookup,lookupCount,LowerBound,LUdecomposition,M0nbar,M2CODE,Macaulay2Doc,makeDirectory,MakeDocumentation,makeDocumentTag,MakeHTML,MakeInfo,MakeLinks,makePackageIndex,MakePDF,makeS2,Manipulator,map,MapExpression,MapleInterface,markedGB,Markov,MarkUpType,match,mathML,Matrix,matrix,MatrixExpression,Matroids,max,maxAllowableThreads,maxExponent,MaximalRank,maxPosition,MaxReductionCount,MCMApproximations,member,memoize,memoizeClear,memoizeValues,MENU,merge,mergePairs,META,method,MethodFunction,MethodFunctionBinary,MethodFunctionSingle,MethodFunctionWithOptions,methodOptions,methods,midpoint,min,minExponent,mingens,mingle,minimalBetti,MinimalGenerators,MinimalMatrix,minimalPresentation,minimalPresentationMap,minimalPresentationMapInv,MinimalPrimes,minimalPrimes,minimalReduction,Minimize,minimizeFilename,MinimumVersion,minors,minPosition,minPres,minprimes,Minus,minus,Miura,MixedMultiplicity,mkdir,mod,Module,module,ModuleDeformations,modulo,MonodromySolver,Monoid,monoid,MonoidElement,Monomial,MonomialAlgebras,monomialCurveIdeal,MonomialIdeal,monomialIdeal,MonomialIntegerPrograms,MonomialOrbits,MonomialOrder,Monomials,monomials,MonomialSize,monomialSubideal,moveFile,multidegree,multidoc,multigraded,MultigradedBettiTally,MultiGradedRationalMap,multiplicity,MultiplicitySequence,MultiplierIdeals,MultiplierIdealsDim2,MultiprojectiveVarieties,mutable,MutableHashTable,mutableIdentity,MutableList,MutableMatrix,mutableMatrix,NAGtypes,Name,nanosleep,Nauty,NautyGraphs,NCAlgebra,NCLex,needs,needsPackage,Net,net,NetFile,netList,new,newClass,newCoordinateSystem,NewFromMethod,newline,NewMethod,newNetFile,NewOfFromMethod,NewOfMethod,newPackage,newRing,nextkey,nextPrime,nil,NNParser,NoetherianOperators,NoetherNormalization,NonAssociativeProduct,NonminimalComplexes,nonspaceAnalyzer,NoPrint,norm,normalCone,Normaliz,NormalToricVarieties,not,Nothing,notify,notImplemented,NTL,null,nullaryMethods,nullhomotopy,nullParser,nullSpace,Number,number,NumberedVerticalList,numcols,numColumns,numerator,numeric,NumericalAlgebraicGeometry,NumericalCertification,NumericalImplicitization,NumericalLinearAlgebra,NumericalSchubertCalculus,numericInterval,NumericSolutions,numgens,numRows,numrows,odd,oeis,of,ofClass,OL,OldPolyhedra,OldToricVectorBundles,on,OneExpression,OnlineLookup,OO,oo,ooo,oooo,openDatabase,openDatabaseOut,openFiles,openIn,openInOut,openListener,OpenMath,openOut,openOutAppend,operatorAttributes,Option,OptionalComponentsPresent,optionalSignParser,Options,options,OptionTable,optP,or,Order,order,OrderedMonoid,orP,OutputDictionary,Outputs,override,pack,Package,package,PackageCitations,PackageDictionary,PackageExports,PackageImports,PackageTemplate,packageTemplate,pad,pager,PairLimit,pairs,PairsRemaining,PARA,Parametrization,parent,Parenthesize,Parser,Parsing,part,Partition,partition,partitions,parts,path,pdim,peek,PencilsOfQuadrics,Permanents,permanents,permutations,pfaffians,PHCpack,PhylogeneticTrees,pi,PieriMaps,pivots,PlaneCurveSingularities,plus,poincare,poincareN,Points,polarize,poly,Polyhedra,Polymake,PolynomialRing,Posets,Position,position,positions,PositivityToricBundles,POSIX,Postfix,Power,power,powermod,PRE,Precision,precision,Prefix,prefixDirectory,prefixPath,preimage,prepend,presentation,pretty,primaryComponent,PrimaryDecomposition,primaryDecomposition,PrimaryTag,PrimitiveElement,Print,print,printerr,printingAccuracy,printingLeadLimit,printingPrecision,printingSeparator,printingTimeLimit,printingTrailLimit,printString,printWidth,processID,Product,product,ProductOrder,profile,profileSummary,Program,programPaths,ProgramRun,Proj,Projective,ProjectiveHilbertPolynomial,projectiveHilbertPolynomial,ProjectiveVariety,promote,protect,Prune,prune,PruneComplex,pruningMap,Pseudocode,pseudocode,pseudoRemainder,Pullback,PushForward,pushForward,Python,QQ,QQParser,QRDecomposition,QthPower,Quasidegrees,QuaternaryQuartics,QuillenSuslin,quit,Quotient,quotient,quotientRemainder,QuotientRing,Radical,radical,RadicalCodim1,radicalContainment,RaiseError,random,RandomCanonicalCurves,RandomComplexes,RandomCurves,RandomCurvesOverVerySmallFiniteFields,RandomGenus14Curves,RandomIdeals,randomKRationalPoint,RandomMonomialIdeals,randomMutableMatrix,RandomObjects,RandomPlaneCurves,RandomPoints,RandomSpaceCurves,Range,rank,RationalMaps,RationalPoints,RationalPoints2,ReactionNetworks,read,readDirectory,readlink,readPackage,RealField,RealFP,realPart,realpath,RealQP,RealQP1,RealRoots,RealRR,RealXD,recursionDepth,recursionLimit,Reduce,reducedRowEchelonForm,reduceHilbert,reductionNumber,ReesAlgebra,reesAlgebra,reesAlgebraIdeal,reesIdeal,References,ReflexivePolytopesDB,regex,regexQuote,registerFinalizer,regSeqInIdeal,Regularity,regularity,relations,RelativeCanonicalResolution,relativizeFilename,Reload,remainder,RemakeAllDocumentation,remove,removeDirectory,removeFile,removeLowestDimension,reorganize,replace,RerunExamples,res,reshape,ResidualIntersections,ResLengthThree,Resolution,resolution,ResolutionsOfStanleyReisnerRings,restart,Result,resultant,Resultants,return,returnCode,Reverse,reverse,RevLex,Right,right,Ring,ring,RingElement,RingFamily,ringFromFractions,RingMap,rootPath,roots,rootURI,rotate,round,rowAdd,RowExpression,rowMult,rowPermute,rowRankProfile,rowSwap,RR,RRi,rsort,run,RunDirectory,RunExamples,RunExternalM2,runHooks,runLengthEncode,runProgram,same,saturate,Saturation,scan,scanKeys,scanLines,scanPairs,scanValues,schedule,schreyerOrder,Schubert,Schubert2,SchurComplexes,SchurFunctors,SchurRings,SCRIPT,scriptCommandLine,ScriptedFunctor,SCSCP,searchPath,sec,sech,SectionRing,SeeAlso,seeParsing,SegreClasses,select,selectInSubring,selectVariables,SelfInitializingType,SemidefiniteProgramming,Seminormalization,separate,SeparateExec,separateRegexp,Sequence,sequence,Serialization,serialNumber,Set,set,setEcho,setGroupID,setIOExclusive,setIOSynchronized,setIOUnSynchronized,setRandomSeed,setup,setupEmacs,sheaf,SheafExpression,sheafExt,sheafHom,SheafOfRings,shield,ShimoyamaYokoyama,short,show,showClassStructure,showHtml,showStructure,showTex,showUserStructure,SimpleDoc,simpleDocFrob,SimplicialComplexes,SimplicialDecomposability,SimplicialPosets,SimplifyFractions,sin,singularLocus,sinh,size,size2,SizeLimit,SkewCommutative,SlackIdeals,sleep,SLnEquivariantMatrices,SLPexpressions,SMALL,smithNormalForm,solve,someTerms,Sort,sort,sortColumns,SortStrategy,source,SourceCode,SourceRing,SPACE,SpaceCurves,SPAN,span,SparseMonomialVectorExpression,SparseResultants,SparseVectorExpression,Spec,SpechtModule,SpecialFanoFourfolds,specialFiber,specialFiberIdeal,SpectralSequences,splice,splitWWW,sqrt,SRdeformations,stack,stacksProject,Standard,standardForm,standardPairs,StartWithOneMinor,stashValue,StatePolytope,StatGraphs,status,stderr,stdio,step,StopBeforeComputation,stopIfError,StopWithMinimalGenerators,Strategy,String,STRONG,StronglyStableIdeals,STYLE,Style,style,SUB,sub,SubalgebraBases,sublists,submatrix,submatrixByDegrees,Subnodes,subquotient,SubringLimit,Subscript,subscript,SUBSECTION,subsets,substitute,substring,subtable,Sugarless,Sum,sum,SumOfTwists,SumsOfSquares,SUP,super,SuperLinearAlgebra,Superscript,superscript,support,SVD,SVDComplexes,switch,SwitchingFields,sylvesterMatrix,Symbol,symbol,SymbolBody,symbolBody,SymbolicPowers,symlinkDirectory,symlinkFile,symmetricAlgebra,symmetricAlgebraIdeal,symmetricKernel,SymmetricPolynomials,symmetricPower,synonym,SYNOPSIS,syz,Syzygies,SyzygyLimit,SyzygyMatrix,SyzygyRows,syzygyScheme,TABLE,Table,table,take,Tally,tally,tan,TangentCone,tangentCone,tangentSheaf,tanh,target,Task,taskResult,TateOnProducts,TD,temporaryFileName,tensor,tensorAssociativity,TensorComplexes,terminalParser,terms,TEST,Test,testExample,testHunekeQuestion,TestIdeals,TestInput,tests,TEX,tex,TeXmacs,texMath,Text,TH,then,Thing,ThinSincereQuivers,ThreadedGB,threadVariable,Threshold,throw,Time,time,times,timing,TITLE,TO,to,TO2,toAbsolutePath,toCC,toDividedPowers,toDual,toExternalString,toField,TOH,toList,toLower,top,top,topCoefficients,Topcom,topComponents,topLevelMode,Tor,TorAlgebra,Toric,ToricInvariants,ToricTopology,ToricVectorBundles,toRR,toRRi,toSequence,toString,TotalPairs,toUpper,TR,trace,transpose,TriangularSets,Tries,Trim,trim,Triplets,Tropical,true,Truncate,truncate,truncateOutput,Truncations,try,TSpreadIdeals,TT,tutorial,Type,TypicalValue,typicalValues,UL,ultimate,unbag,uncurry,Undo,undocumented,uniform,uninstallAllPackages,uninstallPackage,Unique,unique,Units,Unmixed,unsequence,unstack,Up,UpdateOnly,UpperTriangular,URL,urlEncode,Usage,use,UseCachedExampleOutput,UseHilbertFunction,UserMode,userSymbols,UseSyzygies,utf8,utf8check,validate,value,values,Variable,VariableBaseName,Variables,Variety,variety,vars,Vasconcelos,Vector,vector,VectorExpression,VectorFields,VectorGraphics,Verbose,Verbosity,Verify,VersalDeformations,versalEmbedding,Version,version,VerticalList,VerticalSpace,viewHelp,VirtualResolutions,VirtualTally,VisibleList,Visualize,wait,WebApp,wedgeProduct,weightRange,Weights,WeylAlgebra,WeylGroups,when,whichGm,while,width,wikipedia,Wrap,wrap,WrapperType,XML,xor,youngest,zero,ZeroExpression,zeta,ZZ,ZZParser}
}
\lstalias{Macaulay2output}{Macaulay2}

\title{State estimation of a moving frequency source from observations at multiple receivers}

\author{
 Michela Mancini \\
  Guggenheim~School~of~Aerospace~Engineering\\
Georgia~Institute~of Technology\\
 Atlanta,~GA,~30332\\
  \texttt{mmancini32@gatech.edu} \\
   \And
 Anton Leykin \\
  School~of~Mathematics\\
  Georgia~Institute~of Technology\\
  Atlanta,~GA,~30332\\
  \And
 John A. Christian\\
Guggenheim~School~of~Aerospace~Engineering\\
Georgia~Institute~of Technology\\
 Atlanta,~GA,~30332\\
}

\begin{document}

\maketitle

\begin{abstract}
The task of position and velocity estimation of a moving transmitter (with either a known or unknown frequency) is a common problem arising in many different application domains. 
Based on the Doppler effect, this work presents a direct solution using only the frequency measured by a multitude of receivers with a known state. 
A natural rewriting of the problem as a system of polynomial equations allows for the use of homotopy continuation to find the global solution without any \emph{a priori} information about the frequency source. 
We show that the data from six or seven receivers is sufficient in case of known or unknown frequency, respectively. 
After a brief development of the mathematics, two simple examples are provided: (1) position and velocity estimation of a vocalizing dolphin emitting an acoustic signal and (2) initial orbit determination of a satellite emitting an electromagnetic signal.
\end{abstract}

\section{INTRODUCTION}
Despite the criticism that C. A. Doppler encountered after publishing his famous principle in 1842 (\cite{Doppler:1842,Nolte:2020}), the Doppler effect has subsequently proved to be an important phenomenon in many different application domains. Indeed, we can find examples of its use in many different disciplines, including medicine (\cite{Kaunitz:2016}), astronomy (\cite{Lambourne:1997}), satellite orbit determination (\cite{Guier:1958}), and many others.
In this work, we consider how the Doppler effect can be exploited to simultaneously determine the position and velocity of an object emitting a signal---acoustic or electromagnetic---with either a known or unknown transmit frequency. 

\begin{figure}
\centering
  \includegraphics[width=1\linewidth]{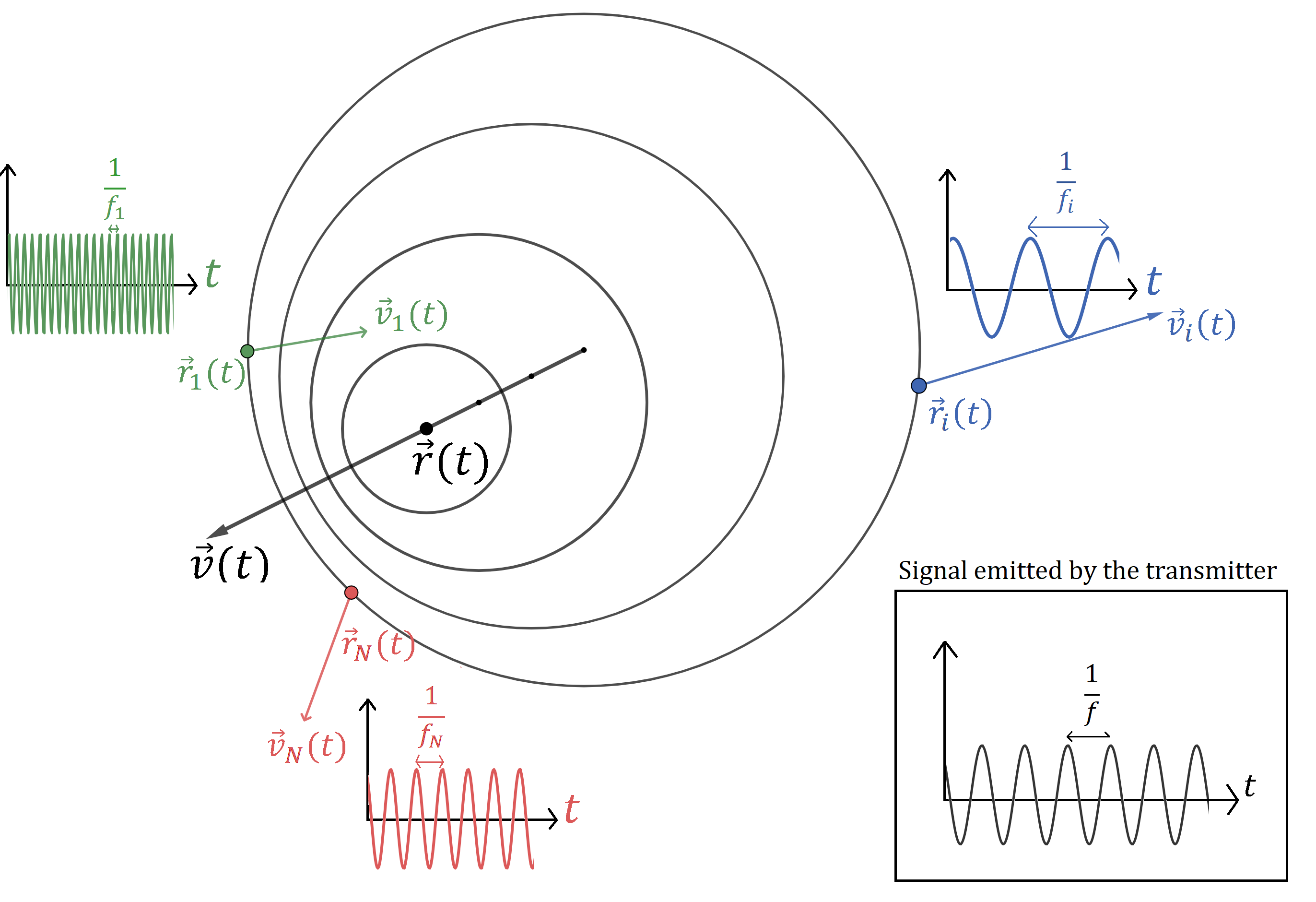}
  \caption{Illustration of the Doppler effect. If the transmitter (black dot) and a receiver have a positive range rate (i.e., are moving away from one another) then the measured frequency will decrease (e.g., blue dot). Conversely, if the transmitter and receiver have a negative range rate (i.e. are approaching one another) then the measured frequecy will increase (e.g., green or red dots).}
  \label{fig:DopplerVisualization}
\end{figure}

Consider a particle with unknown state (position and velocity) that emits a signal. Suppose this signal is subsequently detected by a receiver with a known state (position and velocity). If there is relative motion between the two objects, then the Doppler effect tells us that the frequency measured by the receiver may be different from the frequency at the transmitter. A visualization of this phenomenon is given in Fig. \ref{fig:DopplerVisualization}.
We will show that the classical equation for the Doppler effect may be rearranged to yield a polynomial in the unknown transmitter state. 
If the frequency at the source is known, six receivers suffice to determine the position and velocity of the transmitter up to finitely many possibilities. If it is unknown, seven receivers are needed, and the unknown frequency can be estimated as well. The polynomial system obtained can be solved directly via homotopy continuation, without the need of any initial guess.

The literature presents many circumstances where localizing an unknown transmitter is encountered. Among the marine applications, for example, it can be used to determine the position and velocity of a vocalizing underwater animal (\cite{Gillespie:2020}). It may also be used to localize nodes within an underwater sensor network, which are widely used to collect data for underwater exploration, environmental monitoring, natural disaster prevention (\cite{Felemban:2015}).
Doppler-based state estimation may also be used to track or navigate moving vehicles on the ground, in air, or in space.

We have examples of many previous solutions in all these fields, including for underwater device detection (\cite{Datta:2021,Gong:2020}), drone localization and tracking (\cite{Chan:1990,Famili:2020}), and initial orbit determination (IOD) (\cite{Patton:1960,Christian:2022,Ertl:2023}). However, to the best of our knowledge, most of the past approaches process Doppler measurements within the framework of a sequential filter (e.g., some variant of a Kalman filter) that requires an initial guess of the transmitter state---with it often being unclear how such an \emph{a priori} state estimate might be obtained. 
An attempt to solve the initialization problem was recently studied for the problem of spacecraft IOD (\cite{Ertl:2023}), but the runtimes for the general problem were prohibitive.

In this work we introduce a fast and simple solution to the generic Doppler-only state estimation problem based on homotopy continuation. This allows the analyst to find the global solution in the absence of \emph{a priori} state information, thus providing a means of obtaining the initial guess required by earlier work. We provide {\tt Macaulay2} source code and illustrate this method on simple underwater and space examples.

\section{Polynomial formulation}
Consider a transmitter at the unknown location $\br \in \RR^3 $ moving with the unknown velocity $\bv \in \RR^3 $ and emitting a signal with a frequency $f \in \RR_{>0}$. Assume that the signal is detected by a receiver, placed at the position \(\br_i\), moving with velocity \(\bv_i\). This situation is shown in Fig. \ref{fig:DopplerVisualization}. Due to Doppler effect, the frequency seen by the receiver \(f_i\) will be related to the frequency at the source through the relation
\begin{equation}\label{eq:frequencyshift}
    f_i=\left(1-\frac{\dot{\rho}_i}{c}\right)f
\end{equation}
where \(c\) is the signal's propagation speed in the given medium, and the range rate \(\dot{\rho}_i\) can be obtained taking the derivative of the range \(\rho_i\) between the transmitter and the receiver
\begin{equation}\label{eq:rho}
    \rho_i=\sqrt{\left(\br_i-\br\right)^T\left(\br_i-\br\right)}
\end{equation}
that yields
\begin{equation}\label{eq:rhodot}
    \dot{\rho}_i=\frac{1}{\rho_i}\left(\br_i-\br\right)^T\left(\bv_i-\bv\right)
\end{equation}

A few simple manipulations lead to the polynomial expression we seek. Rearrangement of Eqs.~\eqref{eq:frequencyshift}--\eqref{eq:rhodot} allows us to write the product \(\rho_i\dot{\rho}_i\) as
\begin{equation}\label{eq:main_not_squared}
    \rho_i\dot{\rho}_i = \sqrt{\left(\br_i-\br\right)^T\left(\br_i-\br\right)}\left(\frac{c\left(f-f_i\right)}{f}\right)=\left(\br_i-\br\right)^T\left(\bv_i-\bv\right)
\end{equation}

Squaring the middle and right-hand term of Eq.~\eqref{eq:main_not_squared}, we obtain a polynomial expression relating the known parameters \( (\br_i,\bv_i,f_i) \) to the unknown parameters \( (\br,\bv,f) \):
\begin{equation}\label{eq:main_equation}
    c^2\left(f-f_i\right)^2\left(\br_i-\br\right)^T\left(\br_i-\br\right)-f^2\left[\left(\br_i-\br\right)^T\left(\bv_i-\bv\right)\right]^2=0
\end{equation}
The repetition of Eq.~\eqref{eq:main_equation} for seven receivers produces a polynomial system of seven equations in seven unknowns. If the transmit frequency \(f\) is known, six receivers may be used to solve for the remaining six unknowns. 

Finally, it is important to note the \(c^2\) coefficient in the left-most term of Eq.~\eqref{eq:main_equation}. In applications with high propagation speeds (e.g., speed of light), scaling can have a significant impact on the ability to solve such polynomial systems on a finite precision computer. In such cases, rescaling of the polynomials may provide considerable performance improvements.

\section{Numerical Solution with Homotopy Continuation}
The polynomial system formed by a multitude of frequency observations---each represented by its own instantiation of Eq.~\eqref{eq:main_equation}---can be solved in many different ways.  One family of techniques is homotopy continuation (\cite{Morgan87,Sommese-Wampler-book-05}). In this work, we use an implementation of parameter homotopy provided by the  software package \texttt{NAG4M2} (\cite{Ley11}) in the computer algebra system \texttt{Macaulay2} (\cite{M2}). 

Consider a square system, \(N\) equations in \(N\) unknowns, of polynomials that we aim to solve via homotopy continuation. The procedure begins by selecting a \emph{start system} \(\mathcal{S}\) which belongs to the same family as the \emph{target system} \(\mathcal{F}\) that we desire to solve---which, in this particular case, is a system formed by multiple copies of Eq.~\eqref{eq:main_equation}. We intelligently choose the start system to be one that is easily solvable and then proceed to compute the corresponding solution (e.g., \cite{DuffMonodromy}). Then, beginning from the known solution to the start system \(\mathcal{S}\), homotopy continuation follows a path (in the space of the systems) from start to target while deforming the corresponding solutions---typically using an adaptive predictor-corrector method. Specifically, we track the solutions of the system
\begin{equation}
    \label{eq:HomotopyContinuation1}
    \mathcal{H} = (1-t)\,\mathcal{S} + t\, \mathcal{F} 
\end{equation}
as the continuation parameter \(t\) is varied from 0 to 1.

Construction of the family of polynomial systems is a key factor in this approach. This is critical because the number of solutions of a generic system within the polynomial family corresponds to the number of continuation paths that must be tracked. For instance, many blackbox solvers implement the idea of \emph{total degree homotopy}. Total degree homotopy, an approach which works for any square system, is attractive because it results in a simple start system whose solutions are easy to write down. However, the number of continuation paths equals the total degree in B\'{e}zout's theorem (see chapter 18 of \cite{Harris:1992}), which is often much larger than the number of target solutions we seek. The consequence is that total degree homotopy often results in more paths than are necessary. Moreover, tracking these excess paths (which are divergent) significantly increases the computational cost. Reliance on total degree homotopy likely led to exceptionally long runtimes in a similar approach attempted recently (\cite{Ertl:2023}).

A more nuanced approach is to use the natural parametric family arising in our specific problem. Once the solutions to a {\em generic} system are known we use \emph{parameter homotopy} to deform them into the target solutions. We can choose a random start system in the family---this way the start and target systems are structurally the same differing only in the values of the parameters. 
Depending on the application, the speedup can be considerable. We created and solved the start system using an approach based on monodromy (\cite{DuffMonodromy}). Once the solutions to the start system are known, they can be used solve any other instance of the problem in the family. A short example code that shows the setup of the polynomial system and the calls necessary to perform the tracking in \texttt{Macaulay2} is provided in the Appendix. 

We now apply these ideas to two different scenarios. The first is the generic case, where both the source and receivers are permitted to move. The second is a simplified special case where the source is moving, but the receivers are stationary.

\subsection{General case: moving receivers}
Consider the polynomial system given by Eq.~\eqref{eq:main_equation} specialized for seven different receivers with an unknown transmit frequency \(f\). According to Bézout's theorem, we can expect up to \(6^7=279,936\) solutions in the form of complex-valued 7-tuples \(\left(\br,\,\bv,\,f\right)\). However, the actual number of solutions that we found with techniques based on Gr\"obner basis (\cite{Cox:2015}) is only \(672\). Using the total degree homotopy, almost all the paths will diverge and many hours (on contemporary laptop computers) are required to converge to the 672 solutions. Instead, parameter homotopy allows us to track only \(672\) solutions and solve the desired polynomial system more quickly. We can repeat a similar reasoning for the case of known frequency \(f\), where the total degree homotopy requires the tracking of \(4^6=4,096\) solutions as 6-tuples \(\left(\br,\,\bv\right)\) , while the parameter homotopy allows us to track the actual number of 128 solutions. 

We have experimentally demonstrated that most of the 672 (for unknown \(f\)) or 128 (for known \(f\)) complex solutions can be eliminated without using observations from any additional receivers. In fact, even if they \textit{all} satisfy the polynomial system obtained with Eq.~\eqref{eq:main_equation}, not all satisfy the original (i.e., non-squared) relation from Eq.~\eqref{eq:main_not_squared}.
However, it is not guaranteed that this will eliminate all but one candidates. If an additional observation is available, we can discard the remaining ``wrong''  solutions by checking the residuals of Eq.~\eqref{eq:main_not_squared} written for the additional receiver. The screening procedure described here was sufficient to uniquely identify the desired solution in all of the experiments we performed.

\subsection{Special case: stationary receivers}
Although the general case encompasses all classes of observer motion (both moving and stationary), the special case of stationary observers (\(\bv_i=\textbf{0}\)) simplifies the polynomial systems and permits a more efficient solution. In particular, setting \(\bv_i=\textbf{0}\) causes some of the terms in the original polynomial vanish. Furthermore, we observe that if \((\br,\,\bv,\,f)\) is a solution to the polynomial system, then the same must be true for \((\br,\,-\bv,\,f)\). Counting the number of solutions for this special case, we found that the number of solutions of the polynomial system with stationary observers is 296 for an unknown transmit frequency \(f\). However, thanks to the symmetry of the problem, we can discard half of them and track only 148 paths. Even better, if the frequency at the source is known, then we only need to track 24 paths. This reduction in the number of paths increases the speed of the solver. 

\section{Experiments}
The efficacy of the proposed solution is demonstrated through two illustrative examples.

\subsection{Example 1: Underwater localization with stationary receivers}
The first illustrative experiment is the simultaneous position and velocity estimation of a vocalizing \emph{tursiops truncatus} (the common bottlenose dolphin) swimming within a volume equipped with a network of acoustic sensors (e.g., hydrophones). We assume that the receivers are stationary, with well-known positions from prior surveys. Recognizing that these dolphins have typical swim speeds of 1-4 m/s (\cite{Bernd:1979}) and their whistles are around 7.4-17.2 kHz (\cite{Oswald:2003}), we simulate a dolphin with the following state that is unknown \emph{a priori} to the sensor network:
\[
\br=\begin{bmatrix}
    -5.23\\ 5.28 \\ -15.00
\end{bmatrix} m, \qquad \bv=\begin{bmatrix}
    1.38\\
    1.53\\ 0.22
\end{bmatrix} m/s, \qquad f=15\,kHz
\]

When provided perfect measurements (i.e., no noise), we recover the dynamical state (position and velocity) of the source to machine precision in less than 0.5 seconds for unknown $f$ ($N=7$ observations) and in less than 0.1 seconds for known $f$ ($N=6$ observations). Due to the simplicity of the problem's mathematical structure, an approach similar to that of \cite{hruby2022learning} may lead to an exceptionally fast solver with execution time less than a millisecond per problem instance. Implementation of such a speed-up was not pursued further here as the current sub-second runtimes were fast enough for our proof-of-concept purposes.

We also conducted a Monte Carlo analysis to study the affect of imperfect receiver state knowledge and noisy frequency measurements. Figure~\ref{fig:USN} shows the results from 1,000 numerical experiments for both a known ($N=6$) and unknown ($N=7$) whistle frequency. In each experiment, we apply additive, zero-mean Gaussian noise to each of the observers' position (\(\sigma_r = 1.5 \,cm\)) and to each of the observed frequencies (\(\sigma_f=0.1\,Hz\)).

\subsection{Example 2: Orbit determination with moving receivers}
The second experiment consists of the initial orbit determination (IOD) of a satellite emitting an electromagnetic signal. One receiver is placed aboard a satellite in a circular mid-Earth orbit (MEO) and all of the remaining receivers are placed on the Earth's surface (and rotate with it). We require one orbiting receiver since the state observability is poor when all the receivers are on the Earth's surface. 
The orbital elements (\cite{Bate:2020}) for the transmitting spacecraft are assumed to be
\[
a=12,000\,km,\qquad e=0.1,\qquad i=20^{\circ},\qquad \Omega=200^{\circ},\qquad \omega=20^{\circ},\qquad \nu=0^{\circ}
\]
where $a$ is the semi-major axis, $e$ is the eccentricity, $\Omega$ is the right ascension of the ascending node, $i$ is the inclination, $\omega$ is the argument of the periapsis and $\nu$ is the true anomaly at epoch.
These orbital elements may be converted to three-dimensional position and velocity using standard equations found in most modern astrodynamics textbooks (\cite{Bate:2020,Vallado:2007,Curtis:2020}). With perfect (noise free) receiver state knowledge and frequency measurements, we are able to recover the 6-Degree-of-Freedom dynamical state (and, hence, the six orbital elements) to nearly machine precision.

\begin{figure}[ht!]
\centering
\begin{minipage}[t]{0.45\textwidth}
  \centering
  \includegraphics[width=1\linewidth]{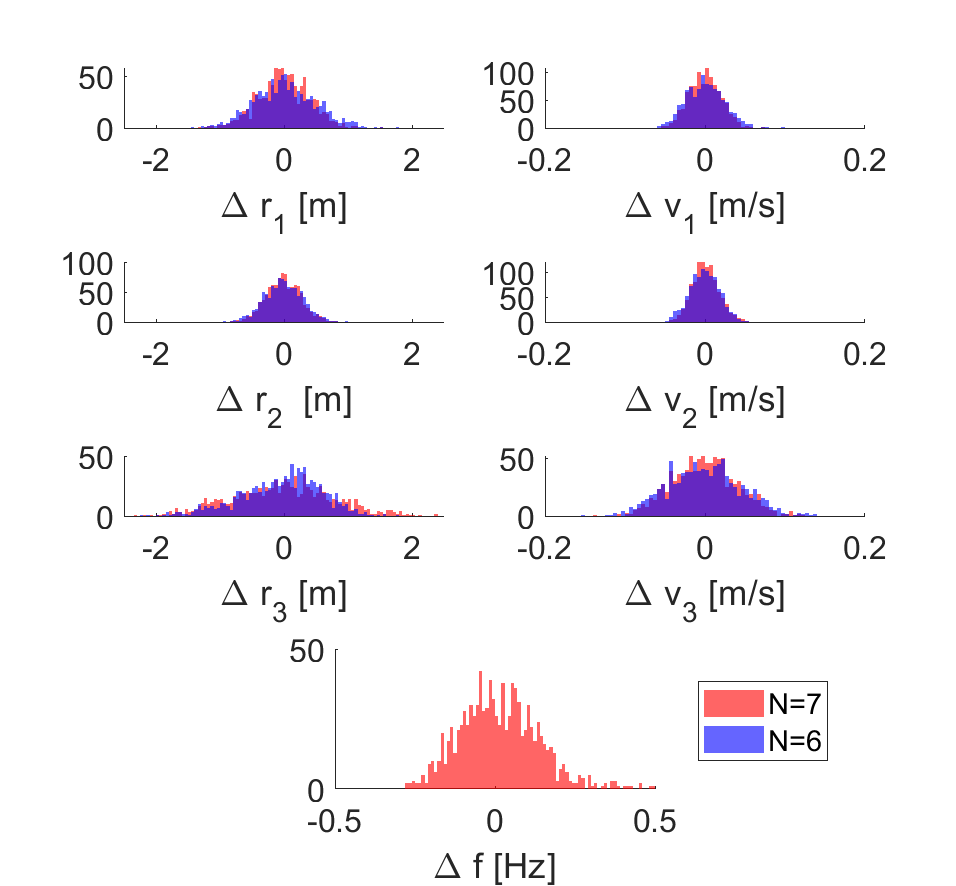}
  \caption{State errors in position (left column) and velocity (right column) for example of a vocalizing dolphin (underwater application). Frequency estimate errors are shown at the bottom for the case of $N=7$ receivers. Results are from a 1,000 trial Monte Carlo analysis.}
  \label{fig:USN}
\end{minipage} \hfil
\begin{minipage}[t]{0.45\textwidth}
  \centering
  \includegraphics[width=1\linewidth]{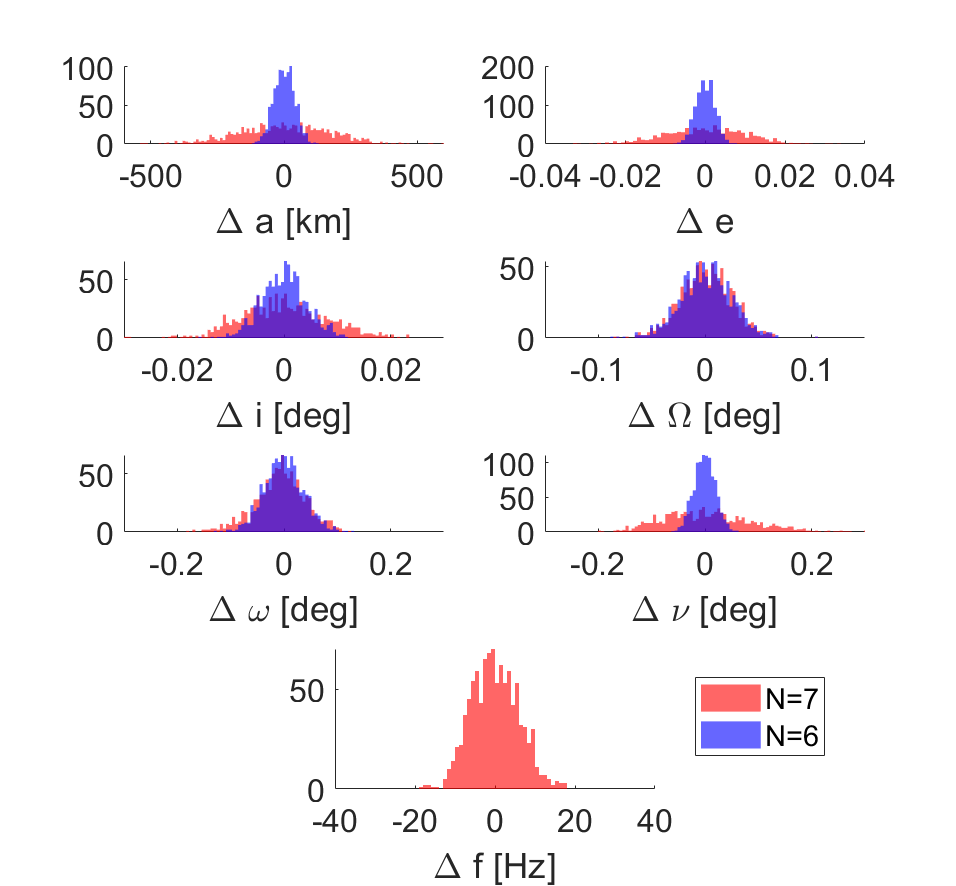}
  \caption{Orbital element errors for example of initial orbit determination of a transmitting satellite (space application). Frequency estimate errors are shown at the bottom for the case of $N=7$ receivers. Results are from a 1,000 trial Monte Carlo analysis. }
  \label{fig:IOD}
\end{minipage}
\end{figure}

As with the underwater case, we may study how performance is affected by imperfect knowledge of the receiver states and errors in the measured frequency. Results from 1,000 numerical experiments are shown in Fig.~\ref{fig:IOD} for both a known ($N=6$) and unknown ($N=7$) transmit frequency. In each experiment, we apply additive, zero-mean Gaussian noise to each of the observers' position and velocity and to the observed frequencies (\(\sigma_f = 0.5 \, Hz\)). Specifically, we used \(\sigma_r = 5 \, cm\) and \(\sigma_v = 1 \, mm/s\) for the receivers on the Earth's surface, and \(\sigma_r = 1 \, m\) and \(\sigma_v = 2 \, cm/s\) for the receiver on orbit.

\section{Conclusions}
The position and velocity of an object emitting an acoustic or electromagnetic signal may be inferred from only the apparent frequency as measured by a multitude of known receivers. This state estimation problem may be written as a system of polynomials, which permits an efficient solution via homotopy continuation. The frequency at the source may be either known or unknown. We highlight the importance of using parameter homotopy, instead of the total degree homotopy, and provide a strategy to filter out multiple solutions. 

Two example applications show the feasibility of this approach, both in the perfect (noise-free) and imperfect (noisy) scenarios. For the imperfect case, we model errors in the receiver states (both position and velocity) and noise is applied to the measured frequency. The proposed method provided a noisy estimate of the transmitter state consistent with the amount of noise in the measurements, both in the case of known and unknown frequency. 

\appendix
\section{Example Code}\label{appendixA}
This appendix provides example code using \texttt{Macaulay2}, which is a popular software system for solving polynomial systems and other problems in algebraic geometry. For more information on \texttt{Macaulay2}, see \cite{M2}. 

The following code example shows how to build a polynomial system based on multiple instantiations of Eq.~\eqref{eq:main_equation}. 

\begin{lstlisting}[language=Macaulay2]
-------
-- import the necessary packages
-------
needsPackage "NumericalAlgebraicGeometry"
needsPackage "MonodromySolver"

-------
-- declare unknowns (corresponding to attributes of freq. source)
-------
variableMatrix = gateMatrix{toList vars(x_1..x_3,u_1..u_3,f)}

-- unknown source position (3D vector)
rr = matrix{{x_1},{x_2},{x_3}}      

-- unknown source velocity (3D vector)
vv = matrix{{u_1},{u_2},{u_3}}      

-- unknown source transmit frequency
ff = matrix{{f}}                    

-------
-- declare the parameters of the system
-------

-- fix the number of receivers for the unknown frequency case 
nObs=7;   

-- receivers' positions (3 x nObs matrix)
posObservers = gateMatrix apply(3, i -> apply(nObs, j -> posi_(i,j)))  

-- receivers' velocities (3 x nObs matrix)
velObservers = gateMatrix apply(3, i -> apply(nObs,j -> veli_(i,j))) 

-- measured frquencies (nObs-D vector)
freqObservers = gateMatrix apply(nObs, i -> {freqi_(i)})    

-- speed of the signal
ci = gateMatrix{{signalSpeed_(i)}}                              

-- variable of all the parameters of the system (1 x 50 vector)
parameterMatrix = gateMatrix{flatten entries(posObservers^{0,1,2})} | gateMatrix{flatten entries(velObservers^{0,1,2}) | flatten entries freqObservers | flatten entries ci}

-------
-- build the system of equations
-------

fsquared = gateMatrix(ff*ff);
csquared = gateMatrix(ci*ci);

DPconstraints = flatten apply(nObs, i-> (
	       rri := gateMatrix posObservers_{i};
	       vvi := gateMatrix velObservers_{i};
	       ffi := gateMatrix{freqObservers#i};
	       dri := rri - rr;
	       dvi := vvi - vv;
	       dfi := ff-ffi;
	       drdv := transpose dri * dvi;
              -- Eq. (5) in the text for the i-th receiver
	       csquared * (transpose dri * dri) * (dfi * dfi)- (drdv * drdv) * fsquared
));

-- create the system described by the given parameters, unknowns and constraints
G = gateSystem(parameterMatrix, variableMatrix, gateMatrix DPconstraints)
\end{lstlisting}

Given the polynomial system \texttt{G} produced by the above code, we may write another script to track the solutions using parameter homotopy continuation.
It is assumed that the parameters describing the start and target system have already been loaded in the \(1 \times (7\,\texttt{nObs} + 1)\) vectors \texttt{startParameters} and \texttt{targetParameters}. The same holds for the solutions of the start system, contained in the matrix \texttt{startSolutions}, of dimension (number of unknowns \(\times\) number of solutions). In practice, we obtain \texttt{startSolutions} using \texttt{MonodromySolver}, a package of \texttt{Macaulay2} based on \cite{DuffMonodromy}.

\begin{lstlisting}[language=Macaulay2]
-------
-- track the solutions
-------

-- declare the continuation parameter t, where t=0 for the start system and t=1 for the target system
declareVariable t;

-- define the homotopy given by the chosen start system and target system: see Eq. (6) for reference
H = sub(gateMatrix G, parameters G, 
        (1-t) * startParameters + t * targetParameters);
Hgate = gateHomotopy(H, vars G, t);

-- track the paths from t=0 to t=1
targetSolutions = trackHomotopy(Hgate, startSolutions);

\end{lstlisting}
\clearpage

\section*{Acknowledgment}
The authors thank Rickey Huang for valuable discussions related to the number of polynomial solutions and the implementation of \texttt{Macaulay2} code. We also thank Tim Duff for assistance with developing the monodromy code used to acquire the solutions to our start system.

\section*{Competing insterests declaration}
Competing interests: the authors declare none.

\bibliographystyle{ieeetr}  
\bibliography{references}

\end{document}